\documentclass[11pt]{article}
\font\smallit=cmti10

\usepackage{amssymb,latexsym,amsmath,epsfig,amsthm, xcolor} 

\usepackage{numprint}
\usepackage{hyperref}

\makeatletter

\renewcommand\section{\@startsection {section}{1}{\z@}
{-30pt \@plus -1ex \@minus -.2ex}
{2.3ex \@plus.2ex}
{\normalfont\normalsize\bfseries\boldmath}}

\renewcommand{\@seccntformat}[1]{\csname the#1\endcsname. }

\makeatother

\newtheorem{theorem}{Theorem}

\newcommand{\N}{\mathbb{N}}
\newcommand{\Z}{\mathbb{Z}}

\newcommand{\DL}{{\mathcal{D}}(L)}
\newcommand{\DLk}{{\mathcal{D}}(L_k)}
\newcommand{\al}{\alpha}
\newcommand{\ao}{\alpha_0}

\newcommand{\zi}{\Z[i]}
\newcommand{\re}{{\text{Re}}}
\newcommand{\im}{{\text{Im}}}
\newcommand{\be}{\begin{equation}}
\newcommand{\ee}{\end{equation}}
\newcommand{\beq}{\begin{eqnarray*}}
\newcommand{\eeq}{\end{eqnarray*}}

\begin{document}

\begin{center}
\uppercase{\bf Extending a problem of Pillai to Gaussian lines}

\vskip 20pt
{\bf Elsa Magness
}\\
{\smallit Department of Mathematics, Seattle University, Seattle, WA 98122, USA}\\
{\tt elsa.magness@gmail.com}\\
\vskip 10pt
{\bf Brian Nugent
}\\
{\smallit Department of Mathematics, Seattle University, Seattle, WA 98122, USA}\\
{\tt bnugent@uw.edu}\\ 
\vskip 10pt
{\bf Leanne Robertson}\\
{\smallit Department of Mathematics, Seattle University, Seattle, WA 98122, USA}\\
{\tt robertle@seattleu.edu}\\ 
\end{center}
\vskip 20pt

\centerline{\bf Abstract}

 Let $L$ be a primitive Gaussian line, that is, a line in the complex plane that contains two, and hence infinitely many, coprime Gaussian integers. 
We prove that there exists an integer $G_L$ such that for every integer $n\geq G_L$ there are infinitely many sequences of $n$ consecutive Gaussian integers on $L$ with the property that none of the Gaussian integers in the sequence is coprime to all  the others.  We also investigate the smallest integer $g_L$ such that $L$ contains a  sequence of $g_L$ consecutive Gaussian integers with this property.  We show that $g_L\neq G_L$ in general.  Also, $g_L\geq 7$ for every Gaussian line  $L$, and we give necessary and sufficient conditions for $g_L=7$ and describe infinitely many Gaussian lines with $g_L\geq 260,000$.  We conjecture that both $g_L$ and $G_L$ can be arbitrarily large. Our results extend a well-known problem of Pillai from the rational integers to the Gaussian integers.

\section{Introduction}\label{intro}

In 1940, Pillai \cite{pillai1} proved that every  sequence  of
16 or fewer  consecutive rational integers contains at least one integer  that
is coprime  to  all  the  others. He also conjectured that for every $n\geq 17$ there exists a sequence of $n$ consecutive integers that does not have this property, and proved that this is the case for $17\leq n\leq430$.  The following year, Brauer \cite{brauer}  proved that Pillai's conjecture holds for all $n\geq 17$, and shortly afterwards Pillai himself published the same result \cite{pillai2}.  
Different proofs have also been provided by Evans \cite{evans1} and Gassko \cite{gassko}.

This problem of Pillai has been extended in two main directions: First, replacing the coprimality requirement with a requirement that the GCD have certain prescribed values and, second, replacing the sequence of consecutive integers with consecutive terms of some other sequence.  We are concerned with the second type of generalization. For results in the first direction, see the papers of Caro~\cite{caro}, Hajdu and Saradha \cite{hajdu1}, and Saradha and Thangadurai \cite{saradha}.

Many authors have extended Pillai's problem to other sequences of integers.  Let $A=(a_1, a_2, a_3, \ldots)$ be a sequence of integers.
Following the literature, we define $g_A$ to be the smallest positive integer $n\geq 2$ such that there are $n$ consecutive terms $a_{k+1}, a_{k+2}, \ldots, a_{k+n}$ of $A$ with the property that none of them is coprime to all the others. Similarly, we define $G_A$ to be the smallest positive integer such that for all $n\geq G_A$ there exist $n$ consecutive terms of $A$ where none  is coprime to all the others.  Then, we have $g_\N=G_\N=17$, by the results of Pillai and Brauer cited above. Notice that $g_A$ and $G_A$ do no exist for every choice of $A$, for example if $A$ is the sequence of primes in $\N$. Also, $g_A$ can exist without $G_A$ existing.  For instance, if $A$ is the sequence of Lucas numbers, then Hajdu and Szikszai \cite{hajdu4} prove that $g_A=171$ but $G_A$ does not exist. If $G_A$ exists, however, then clearly $g_A$ also exists and $G_A\geq g_A$.

Following 
Hajdu and Szikszai \cite{hajdu2}, we call $A$ a {\em Pillai sequence} if $G_A$ exists.  Thus, $\N$ is a Pillai sequence, but the sequences of primes and of Lucas numbers are not.  Other examples of Pillai sequences include arithmetic progressions~\cite{evans2} and  non-degenerate linear and elliptic divisibility sequences \cite{hajdu2, hajdu3}. Hajdu and Szikszai \cite{hajdu4} also give a simple characterization of which associated Lucas and Lehmer sequences are Pillai sequences and which are not (the characterization depends only on the parities of the coefficients in their defining recurrences).

More recently,  several authors have applied Pillai's problem to  sequences arising from polynomials. Let $f\in\Z[x]$ be a  polynomial of degree greater than one, and consider the integer sequence $A_f=(f(1), f(2), f(3),\ldots)$. To simplify notation, let $g_f=g_{A_f}$ and $G_f=G_{A_f}$.  Harrington and Jones \cite{harrington}  calculate $g_f$ for various families of quadratic polynomials, and conjecture that $g_f$ exists and is smaller than $35$ for every quadratic polynomial $f$ (they do not consider $G_f$). 
Sanna and Szikszai~\cite{sanna} prove the first part of this conjecture as a corollary of their main result that $G_f$ (and hence $g_f$) exists if  $f$ is quadratic or cubic.  Moreover, they show that for every $n\geq G_f$ there are {infinitely} many integers $k\geq 0$ such that none of the consecutive terms $f(k+1), f(k+2), \ldots, f(k+n)$ of $A_f$ is coprime to all the others.
Ford et al. \hskip-.04in \cite[Corollary~3 of Theorem~1]{ford} extend this  result to {\em all} non-constant polynomials as one of the consequences of their powerful  main theorem that provides a new method for obtaining long gaps in sieved sets.  Specifically, they prove that if $f\in\Z[x]$ is a non-constant polynomial, then there exists an integer $G_f\geq 2$ such that for every integer $n\geq G_f$ there are {\em infinitely} many integers   $k\geq 0$ such that none of the $n$ integers  $f(k+1), \ldots, f(k+n)$  is coprime to all the others. 
Thus, the sequence $A_f=(f(1), f(2), f(3),\ldots)$ is a Pillai sequence for all non-constant polynomials $f\in\Z[x]$.

In this paper, we extend Pillai's problem to sequences of consecutive  Gaussian integers along lines in the complex plane. Gethner et al. \hskip-.04in \cite{gethner} coined the term  {\em{Gaussian line}} for a line in the complex plain that contains two, and hence infinitely many, Gaussian integers.  We call such a line {\em{primitive}}  it it contains two coprime Gaussian integers
(coprimality is well defined in the ring $\Z[i]$ of Gaussian integers since it is a unique factorization domain).  Many of the properties of the rational integers  on the real line extend to the Gaussian integers on a primitive Gaussian line, including the existence of arbitrarily long sequences of composites, the periodicity of divisibility, and the Chinese remainder theorem (see \cite{gethner, magness}).  Here we show that we can include  being a Pillai sequence to the list of shared properties of the rational integers  on the real line and  the Gaussian integers on a primitive Gaussian line.

Our approach is different from that taken by Ghorpade and Ram \cite{ghorpade} who also extend Pillai's problem to the Gaussian integers, and more generally to  arbitrary integral domains where the notion of coprimality is defined (eg., rings of integers of number fields of class number one).  The main difference is that they do not start with an infinite sequence of Gaussian integers or of elements in some other ring $R$, but instead consider the entire ring at once and look for finite arithmetic progressions anywhere in the ring where none of the terms is coprime to all the others.  For a large class of rings $R$ (including unique factorization domains of characteristic zero), Ghorpade and Ram give a formula for a number $N_R$ with the property that any arithmetic progression in $R$ of at most $N_R$ terms with the first term coprime to the common difference contains a term that is coprime to all the other terms.  They also show that if an additional condition on the ring $R$ is satisfied, then for all $n> N_R$ there is an arithmetic progression of $n$ terms of $R$ with the first term coprime to the common difference such that none of the terms is coprime to all the others.  In the special case where $R$ is the ring $\zi$ of Gaussian integers, their formula gives that $N_{\zi}=6$.  This means that any arithmetic progression of six or fewer Gaussian integers where at least two are coprime contains a term that is coprime to all the others, and if $n>6$ then there exists an arithmetic progression of $n$ terms of $\zi$ such that none of the terms is coprime to all the others. Unlike our case, they do not require that all the arithmetic progressions be on the same line.  In particular,  for $n\geq 17$ they can simply use that there exists such an arithmetic progression in $\N$ (and hence in $\zi$) since $\N$ is a Pillai sequence with $G_\N=17$.

An outline of our paper and main results is as follows.  In Section 2 we give an overview of the properties of Gaussian lines that are needed for our work on Pillai's problem applied to these lines. All of the results in this section, except Theorem~\ref{newchinese}, are proven in \cite{magness} but are stated here 
to make the paper more self-contained.  Theorem \ref{newchinese}  is an extension of Theorem 9 in \cite{magness} and shows that there are infinitely many Gaussian lines that satisfy any finite list of divisibility conditions, as long as those conditions do not violate the periodicity of prime divisibility on the line.
We then turn to the problem of Pillai.  For a primitive Gaussian line  $L$ we define a particular Gaussian integer $\alpha_0$ on $L$ and consider the infinite sequence   ($\alpha_0, \alpha_1, \alpha_2, \ldots $)
 of consecutive Gaussian integers on $L$ beginning with $\alpha_0$ and continuing with increasing real part.
We define
$g_L$ to be the smallest positive integer $n\geq 2$ such $L$ contains a sequence of $n$ consecutive Gaussian integers $\alpha_{k+1}, \alpha_{k+2}, \ldots \alpha_{k+n}$  where none of the terms of the sequence is coprime to all the others. Similarly, we define $G_L$
to be the smallest positive integer such that for all
$n\geq G_L$ there is a sequence of $n$ consecutive Gaussian integers $\alpha_{k+1}, \alpha_{k+2}, \ldots \alpha_{k+n}$ on $L$ where none of the terms is coprime to all the others. 
Then $g_L$ and $G_L$ exist in the special case where $L$ is the real line, and in this case we have $g_L=g_\N=G_\N=G_L=17$. In Section 3, we prove that $G_L$ (and hence $g_L$) exists for every primitive Gaussian line $L$. Thus, all Gaussian lines are Pillai sequences.  In Section 4, we turn our attention to $g_L$.  We show that $g_L\neq G_L$ in general. It follows from the result of Ghorpade and Ram~\cite{ghorpade} mentioned above that  $g_L\geq 7$ for all $L$.  We prove this in a different way and classify those Gaussian lines for which $g_L=7$. Considering lines with $g_L=7$, we show 
that for any $B\in\N$ there are infinitely many Gaussian lines $L$ with $g_L=7$ and such that for all $7\leq n\leq B$ there are infinitely many  integers   $k\geq 0$ such that none of the $n$ consecutive Gaussian integers $\alpha_{k+1}, \alpha_{k+2}, \ldots \alpha_{k+n}$ on $L$   is coprime to all the others (thus either $G_L=7$ or $G_L>B$).
Considering large values of $g_L$, we describe infinitely many Gaussian lines for which $g_L\geq 260,000$.  It seems plausible that $g_L$ (and hence $G_L$) can be arbitrarily large.   We leave the reader with this and several other open problems.
 
 \bigskip

\section{Gaussian Lines}\label{gaussian}

In this section we provide the notation  and background about Gaussian lines that  we use for extending Pillai's problem.  All of the results, except Theorem~\ref{newchinese}, are proven in~\cite{magness} but are included here for the convenience of the reader. Throughout, denote the {\em norm} of a Gaussian integer $x+iy$ by $N(x+iy)=x^2+y^2$.

Let $L$ be a Gaussian line.  We distinguish two Gaussian integers, $\ao=a+bi$   and $\delta=c+di$, that define $L$ and provide a natural ordering of the Gaussian integers on $L$ for studying Pillai's problem.
These numbers are defined as follows. Let $\ao$ be the Gaussian integer on $L$ of minimum norm, and if there are two such integers, let $\ao $ be the one with the larger real part.  If $L$ is vertical, then take $\delta=i$.  Otherwise, let $\alpha_1$ be the Gaussian integer on $L$ closest to $\ao$ (so  $N(\alpha_1-\ao)$ is minimal) and with $\re(\alpha_1)>\re(\ao)$.  Then take $\delta=\alpha_1-\ao$.  Thus, $\alpha_0$ is on the line $L$, but  $\delta$ is not, provided $\ao\neq 0$. 
In Lemma 1 in \cite{magness}, we show that $c$ and $d$  are  coprime, $c\geq 0$, and  the Gaussian integers  on $L$ are exactly the Gaussian integers $\alpha_k$ given by
\begin{equation}\label{an}\alpha_k=\ao+\delta k,\  k\in\Z.\end{equation}
Moreover, $L$ is primitive if and only if $\ao$ and $\delta$ are coprime in $\Z[i]$.
In this paper, we study Pillai's problem for the sequence 
$(\alpha_0, \alpha_1, \alpha_2,\ldots )$ 
of consecutive Gaussian integers on a primitive Gaussian line $L$.

From now on, assume that $L$ is a primitive Gaussian line with $\ao$ and $\delta$ defined as above.
We also define a rational integer $\Delta$ associated to $L$ by 
\begin{equation}\label{D}\Delta=ad-bc\in\Z.\end{equation}  
Then,  $\Delta=0$ if and only if $L$ is the real line $\im(z)=0$ or the imaginary line $\re(z)=0$, which holds if and only if $\alpha_0=0$ (see \cite[Lemma 2]{magness}).

Together, $\Delta$ and $\delta$ enable us to easily classify those Gaussian primes that divide some Gaussian integer on $L$.  
Define the {\em divisor set of $L$}, denoted $\DL$, to be the set of Gaussian integers that divide some Gaussian integer on $L$ (so a Gaussian integer in $\DL$ does not necessarily lie on the line $L$, but simply divides some Gaussian integer that lies on $L$). Theorem 11 in \cite{magness} gives a complete characterization of $\DL$ based only on the values of $\Delta$ and $\delta$. 
In this paper we are only concerned with knowing which Gaussian primes occur in $\DL$ since we are focused on the  coprimality of elements on the line.   We will use the following simple test for which rational primes and non-rational Gaussian primes occur in the divisor set $\DL$ of a primitive Gaussian line $L$. Recall that a rational prime $p$ is a Gaussian prime if and only if $p \equiv 3 \pmod 4$.

\begin{theorem}[{\cite[Theorems 4 and 5]{magness}}] \label{DL} Let  $L$ be a primitive Gaussian line and $\DL$ be its divisor set.
\begin{enumerate}
\item[(a)]  If $p\in\Z$ is a rational prime, then $p\in\DL$  if and only if $p$  divides~$\Delta$.
\item[(b)] If $\pi\in\zi$ is a non-rational Gaussian prime, then $\pi\in\DL$ if and only if $\pi$ does not divide $\delta$. 
\end{enumerate}\end{theorem}

\noindent{\em Remark 1.} If $L$ is not the real or imaginary line (so $\Delta\neq 0$), then it follows from Theorem \ref{DL} that there are only finitely many rational primes $p$ that divide some Gaussian integer on $L$. It also follows  that the divisor set $\DL$ contains at least one Gaussian prime lying over $p$ for every rational prime $p\equiv 1\pmod 4$ since 
$\delta=c+di$ is not divisible by any rational integers (otherwise $c$ and $d$ would not be coprime as required).  Specifically, for every rational prime $p\equiv 1\pmod 4$, if $p$ does not divide $N(\delta)$, then $\DL$ contains exactly two non-associate Gaussian primes lying over $p$, and if $p$ divides $N(\delta)$ then it contains exactly one such prime.

\bigskip

The set of Gaussian integers on a primitive Gaussian line shares several properties with the set of rational integers on the real line (see \cite{gethner, magness}).  The two main shared properties important here are the periodicity of divisibility and the Chinese remainder theorem.  Before stating these results for Gaussian lines, we need to define a function $\nu:\Z[i]\rightarrow \Z$ by 
\begin{equation}\label{morm}
\nu(x+iy) = \frac{N(x+iy)}{{\rm gcd}(x,y)}.
\end{equation}
 Notice that $\nu(r)=r$ for all $r\in \Z$. The function $\nu$  is useful  because if $\beta\in\Z[i]$ then the smallest positive rational integer divisible by $\beta$ is $\nu(\beta)$, and furthermore, $\nu(\beta)$ divides every rational integer that is divisible by $\beta$. In particular, if $p$ is a rational prime, then $\nu(\beta)=p$ if and only if $\beta$ is a Gaussian prime lying over $p$. Furthermore, our next theorem shows that divisibility on $L$ by $\beta$ is periodic with period $\nu(\beta)$ for all $\beta\in\DL$.  Thus, the periodicity of divisibility of rational integers on the real line extends to Gaussian lines.

\begin{theorem}[{\cite[Theorem 3]{magness}}]\label{periodicity}
Suppose $\beta\in\Z[i]$ divides some Gaussian integer $\alpha_t$ on $L$.  
Then  $\beta$ divides $\alpha_{k}$ if and only if $k \equiv t \pmod {\nu(\beta)}$.
\end{theorem}

The next theorem extends the Chinese remainder theorem to Gaussian lines.

\begin{theorem}[{\cite[Theorem 8]{magness}}] \label{crtGL}  Let $L$ be a primitive Gaussian line, and suppose
$\mu_1, \mu_2,\ldots, \mu_k$ are Gaussian integers in the divisor set $\DL$ of $L$ such  that $\nu(\mu_1),\nu(\mu_2), \ldots, \nu(\mu_k)$ are pairwise coprime. Let $b_1, b_2,\ldots, b_k\in\Z$. Then there is a unique rational integer $t$ modulo the product $\nu(\mu_1)\nu(\mu_2) \cdots \nu(\mu_k)$ such that
$$ \mu_1\mid\alpha_{t+b_1}, \ \mu_2\mid\alpha_{t+b_2}, \ \ldots, \ \mu_k\mid\alpha_{t+b_k}.$$
\end{theorem}

We will also need a new theorem about Gaussian lines that strengthens Theorem~9 in~\cite{magness}.  This theorem shows that if you want a Gaussian line  with the property that certain Gaussian integers divide specified elements on the line and  certain Gaussian  primes do not divide any elements on the line, then you are in luck as long as a coprimality condition is satisfied.  Namely, not only is there a line that satisfies your desired conditions, but there are infinitely many.

\begin{theorem} \label{newchinese}
Let $\mu_1, \mu_2, \ldots, \mu_k$ be Gaussian integers, 
$b_1, b_2, \ldots, b_k$ be rational integers  (not necessarily distinct), 
$p_1,p_2,\ldots, p_t$ be rational primes that are congruent to 3 modulo 4 (so also Gaussian primes), and 
$\pi_1,\pi_2,\ldots, \pi_\ell$ be non-rational Gaussian primes. 
Suppose that all the $\mu_j$, $p_n$, $\pi_m$ are pairwise coprime ($1\leq j\leq k$, $1 \leq n\leq t$, $1 \leq m\leq \ell$).
Then there are infinitely many primitive Gaussian lines  $L$ that satisfy the following 
 divisibility 
 properties:
 \begin{quote} \begin{enumerate}
 \item[(a)]  $\mu_1, \mu_2, \ldots, \mu_k\in\DL$;
 \item[(b)] 
 $p_1,p_2,\ldots, p_t,\pi_1,\pi_2,\ldots, \pi_\ell\not\in\DL$;
  \item[(c)] $\mu_j $ divides the Gaussian integer $ \alpha_{b_j}$ \hskip-.03in on $L$ for $1 \leq j\leq k$.
 \end{enumerate}
\end{quote}
\end{theorem}
\smallskip

\begin{proof}  
To show there are infinitely many primitive Gaussian lines $L$ that satisfy Properties (a)--(c) in the theorem, we show that there are infinitely many Gaussian integers
 $\alpha_0=a+bi$ and $\delta=c+di$ that satisfy the following six properties:
\begin{quote} \begin{enumerate}
 \item $N(\ao+n\delta)>N(\ao)$ for all $n\neq 0$, $n\in\Z$;
 \item gcd$(c,d)=1$ and $c\geq 0$;
  \item $\alpha_0$ and $\delta$ are coprime over $\zi$;
  \item $\mu_j $ divides  $ \alpha_{b_j}\hskip-.035in=\ao+b_j\delta$  for  $1 \leq j\leq k$;
   \item $p_n $ does not divide  $\Delta$  for $1 \leq n\leq t$;
    \item $\pi_m $ divides  $\delta$  for $1 \leq m\leq \ell$.
 \end{enumerate}
\end{quote}
This is sufficient since it follows from Properties 1--3 that $\alpha_0$ and $\delta$ define a primitive Gaussian line and from Properties 4--6 that this line satisfies Properties (a)--(c) stated in the theorem.

We first choose $\ao$. For $1 \leq j\leq k$, let $\gamma_j\in\Z[i]$ be a greatest common divisor of $\mu_j$ and $b_j$ (so each $\gamma_j$ is uniquely defined up to multiplication by a unit in $\zi$). Let $\lambda$ be a Gaussian integer that is coprime to $\mu_j, b_j, p_n, \pi_m$ for $1\leq j\leq k$, $1 \leq n\leq t$, $1 \leq m\leq \ell$.
Define $\ao$ by 
$$\alpha_0=\lambda\prod_{j=1}^k \gamma_j\in\zi.$$
There are infinitely many possibilities for $\alpha_0$ since there are infinitely many choices for $\lambda$. For each $\alpha_0$, we show there are infinitely many $\delta$ 
for which Properties 1--6 are satisfied for $\ao$ and $\delta$. To do this, we find congruences that are sufficient for Properties 3--6 to hold, then construct infinitely many solutions $\delta$ to these congruences for which Properties 1 and 2 also hold. 

We begin with Property 4. For $1 \leq j\leq k$, the Gaussian  integers $b_j/\gamma_j$ and $\mu_j/\gamma_j$ are coprime, so there exists  $\eta_j\in\Z[i]$ such that
\be\label{eta}1\equiv\eta_j\left(\frac{b_j}{\gamma_j}\right)\pmod{\frac{\mu_j}{\gamma_j}}.\ee
We claim that  Property 4 is satisfied if $\delta$ is a solution to the following system of $k$ congruences:
\be\label{P4}x\equiv -\eta_j\left(\frac{\ao}{\gamma_j}\right)\pmod{\frac{\mu_j}{\gamma_j}}, \ \ \ 1 \leq j\leq k.\ee
The moduli in this system are pairwise coprime, so this system has a solution by the Chinese remainder theorem for the Gaussian integers. If 
 $\delta$ is such  a solution, then
multiplying the congruence in (\ref{eta}) by the Gaussian integer $\alpha_0/\gamma_j$ gives
$$\frac{\alpha_0}{\gamma_j}\equiv -\delta\left(\frac{b_j}{\gamma_j}\right)\pmod{\frac{\mu_j}{\gamma_j}}.$$
By multiplying  by $\gamma_j$  it follows that $\alpha_0\equiv -b_j\delta\pmod{\mu_j}$.
Thus,  $\alpha_j=\alpha_0+b_j\delta$ is divisible by $\mu_j$
for all $1\leq j\leq k$, and Property 4 is satisfied.

 For Property 3, we want a solution to (\ref{P4})  that is also 
  coprime to $\alpha_0$.  Let $\beta$ be the product of all the Gaussian primes that divide $\alpha_0$ and are coprime to $\mu_j/\gamma_j$, $1\leq j\leq k$, and  let $\beta=1$ if no such Gaussian primes exist.  
 Then for Property 3, we want $\delta$ to be coprime to both $\beta$ and 
 $\mu_j/\gamma_j$, $1\leq j\leq k$. But if $\delta$ satisfies (\ref{P4}) then it is automatically coprime to $\mu_j/\gamma_j$, $1\leq j\leq k$, since each $\ao/\gamma_j$ is coprime to $\mu_j/\gamma_j$. Thus, to ensure that $\delta$ is also coprime to $\beta$ we include the sufficient requirement that $\delta\equiv 1 \pmod\beta$.
 
 To ensure that Property 5 is satisfied, we require that 
 \be\label{P5}\delta\equiv i\alpha_0\equiv -b+ai\pmod P,\ee 
 where $P=p_1p_2\cdots p_t$. To see that this is sufficient, first notice that if (\ref{P5}) is satisfied then $\Delta\equiv a^2+b^2=N(\alpha_0)\pmod P$. Thus, if $p_n$ divides $\Delta$ for some $1\leq n\leq t$,
 then $p_n$ divides $N(\alpha_0)$ and so $p_n$ divides $\alpha_0$ since the prime $p_n$ is inert in $\zi$.  This is not possible since $\alpha_0$ is coprime to $P$ by construction.  Thus, $p_n $ does not divide~$\Delta$  for all $1 \leq n\leq t$.

 Finally,  notice that Property 6 is equivalent to $\delta\equiv 0\pmod{\pi_1\pi_2\cdots \pi_\ell}$ since $\pi_1,\pi_2,\ldots, \pi_\ell$ are pairwise coprime.
 Thus,  for $\delta$ to satisfy Properties 3--6,  it is sufficient that $\delta$  be a solution to the following system of $k+3$ congruences:
\begin{align*}\label{P3456}
x &\equiv -\left(\frac{\ao}{\gamma_j}\right)\kappa_j^{-1} \pmod{\frac{\mu_j}{\gamma_j}}, \ \ 1 \leq j\leq k,\\
x & \equiv 1\pmod\beta,\\
x&\equiv i\alpha_0\pmod P, \  \ {\text{and}} \\
x&\equiv 0\pmod{\pi_1\pi_2\cdots \pi_\ell}.
\end{align*}
The moduli in these $k+3$ congruences are pairwise coprime;  let $M$ be their product.  The system has a  unique solution 
 modulo $M$ by the Chinese remainder theorem for the Gaussian integers. Let $\tau\equiv r+si\pmod M$ denote this unique solution.  
 
 It remains to construct $\delta=c+di$ that satisfies  Properties 1 and 2, and such that $\delta\equiv\tau\pmod{M}$, so that Properties 3--6 hold as well.  This follows exactly as in the last part of the proof of Theorem 9 in \cite{magness} except that the modulus $\beta\omega_1\omega_2\cdots\omega_k$ used there is now replaced by $M$ throughout.  This completes the proof of Theorem~\ref{newchinese}.
\end{proof}

\section{Gaussian Lines are Pillai Sequences}

In this section we prove that $G_L$, and hence $g_L$, exists for every primitive Gaussian line $L$. Thus, all Gaussian lines are Pillai sequences, since non-primitive lines trivially have $G_L=g_L=2$. Our proof relies on  the main theorem by Ford, Konyagin, Maynard, Pomerance, and Tao~\cite{ford} about long gaps in sieved sets.

\begin{theorem}\label{GL} Let $L$ be a primitive Gaussian line.  Then there exists an integer $G_L$ such that for all integers $n\geq G_L$ there are infinitely integers $k\geq 0$ with the property that none of the $n$ consecutive Gaussian integers $\alpha_{k+1}, \alpha_{k+2},\ldots, \alpha_{k+n}$ on $L$ is coprime to all the others.

\end{theorem}

\begin{proof}  Let $L$ be a primitive Gaussian line with associated $\alpha_0,\delta\in\Z[i]$. Then $\alpha_0$ and $\delta$ are coprime and the Gaussian integers on $L$ are precisely the numbers $\alpha_x=\alpha_0+x\delta$, $x\in\Z$. We modify the proof of Corollary 3 of Theorem 1 in \cite{ford} to prove that for all large 
$n$ there are infinitely many integers $k$ such that for each $t\in\{1,\ldots,n\}$ there is a $j\in\{1,\ldots,n\}$ with $t\neq j$ such that $\alpha_{k+t}$ and $\alpha_{k+j}$
have a common Gaussian prime divisor $\pi\in\Z[i]$ with $N(\pi)>2$.  Thus, the sequence $\alpha_{k+1}, \alpha_{k+2},\ldots, \alpha_{k+n}$ of $n$ consecutive Gaussian integers on $L$ does not contain an element that is coprime to all the others.

The norm of an arbitrary Gaussian integer $\alpha_x$ on $L$ 
 can be viewed as a quadratic polynomial $f\in\Z[x]$ as follows:
\beq\label{N} f(x) \ = \ N(\alpha_x)&=&N(\alpha_0+\delta x)\\
&=&N(\delta)x^2+Tr(\alpha_0\overline\delta)x+N(\alpha_0)\\
&=&(c^2+d^2)x^2+2(ac+bd)x+a^2+b^2\nonumber.
\eeq
Let $\Delta$ be defined as in Equation (\ref{D}). Then the  discriminant  of $f$ is equal to $-4\Delta^2$, and so is negative unless $\Delta=0$.  If  
$\Delta=0$, then $L$ is the real or imaginary line and Theorem \ref{GL} holds by Brauer's result \cite{brauer}.  Thus, we may assume $\Delta\neq 0$, in which case we have that $f$ is irreducible over $\Z$.

Following the proof of Corollary 3 in \cite{ford}, we let $I_2$ be the empty set, and for all  odd rational prime $p$, we define 
$$I_p=\{n\in\Z/p\Z:\, f(n)\equiv 0\hskip-.07in\pmod p\}.$$
Then, $\left| I_2\right|=0$ and, by Theorem \ref{DL}, we have the following  for odd primes $p$:
\begin{align*}p\nmid N(\delta),  \ p\equiv 1\hskip-.07in\pmod 4&\Longrightarrow \left| I_p\right|=2;\\
p| N(\delta),  \ p\equiv 1\hskip-.07in\pmod 4&\Longrightarrow\left| I_p\right|=1;\\
p\nmid \Delta,  \ p\equiv 3\hskip-.07in\pmod 4&\Longrightarrow \left| I_p\right|=0;\\
p| \Delta,  \ p\equiv 3\hskip-.07in\pmod 4&\Longrightarrow\left| I_p\right|=1.
\end{align*}
For $x\in\N$, consider the {\it sifted set} 
$S_x$
of all  integers $m$ such that for all primes $p\leq x$, the class $\overline{m}\in\Z/p\Z$ is not contained in $I_p$.  That is, $m\in S_x$ if  and only if the norm $N(\alpha_m)$ of the Gaussian integer $\alpha_m$ on $L$ is not divisible by any odd prime $p\leq x$.  
By Theorem~1 in \cite{ford}, there is a bound $B_L$ such that for all $x\geq B_L$ the set $S_x$ contains a gap of length larger than $m=\lfloor{2x}\rfloor$. That is, there is an integer $k$ such that each of the integers $f(k+1), f(k+2),\ldots, f(k+m)$ has an odd prime factor $p\leq x$. Thus, each of the consecutive Gaussian integer $\alpha_{k+1},\alpha_{k+2}, \ldots,\alpha_{k+m}$ on $L$ has a  Gaussian prime factor $\pi$ that lies over an odd rational prime $p\leq x$.  For $j\in\{1, \ldots m\}$, take a Gaussian prime divisor $\pi$ of $\alpha_{k+j}$ that lies over an odd rational prime $p\leq x$. Then $\pi$ divides every $p$th Gaussian integer on $L$ by Theorem~\ref{periodicity}, so in particular, $\pi$ divides both  $\alpha_{k+j+p}$ and $\alpha_{k+j-p}$. But either $\alpha_{k+j+p}$ or $\alpha_{k+j-p}$ is a term of the sequence $\alpha_{k+1},\ldots,\alpha_{k+m}$ since $p\leq x$ and the sequence  has length $m=\lfloor{2x}\rfloor$.  Thus, $\alpha_{k+j}$ shares a common Gaussian prime divisor with at least one other term of the sequence. Since this is true for all $j\in\{1, \ldots m\}$, it follows that the sequence $\alpha_{k+1}, \alpha_{k+2},\ldots, \alpha_{k+m}$ of $m=\lfloor{2x}\rfloor$ consecutive Gaussian integers on $L$ does not contain an element that is coprime to all the others.  
Moreover, there are in fact infinitely many values of $k$ with this property since, by Theorem \ref{periodicity}, divisibility on $L$ by all Gaussian primes of norm not exceeding $x$ is periodic with period 
$$P(x):=\prod_{\substack{p\leq x \\ I_p\neq 0}}p.$$
Taking $G_L=2B_L$, we have that for $n\geq G_L$ there are infinitely many integers $k$ such that none of the $n$ consecutive Gaussian integers $\alpha_{k+1}, \alpha_{k+2},\ldots, \alpha_{k+n}$ on $L$ is coprime to all the others.
\end{proof}

\noindent{\em Remark 2.} Although our proof of Theorem \ref{GL} involves quadratic polynomials,  Sanna and Szikszai's  \cite{sanna} theorem and proof  that $G_f$ exists for all quadratic polynomials $f\in\Z[x]$ is not sufficient for our proof.  It follows from their work and the ideas above that if 
$L$ is a primitive Gaussian line, then there is an integer $G_L$ such that for all $n\geq G_L$ there are infinitely integers $k\geq 0$ with the property that none of the $n$ rational  integers $N(\alpha_{k+1}), N(\alpha_{k+2}),\ldots, N(\alpha_{k+n})$ on $L$ is coprime to all the others.  This is not sufficient for our proof because if $p\equiv 1 \pmod 4$ is a common rational prime divisor of  $N(\alpha_{k+t})$ and $N(\alpha_{k+j})$ then
$\alpha_{k+t}$ and $\alpha_{k+j}$ could still be coprime in $\Z[i]$ because they could be divisible by non-associate primes lying over $p$. By the main theorem by Ford et al., we can assume that they are a distance $p$ apart and this guarantees that they are divisible by the same prime lying over $p$ as needed.

\section{More on Consecutive Integers along Gaussian Lines}

Let $L$ be a primitive Gaussian line and  $g_L$ be the smallest integer such that $L$ contains a  sequence of $g_L$ consecutive Gaussian integers where none is coprime to all the others. Then $g_L$ exists since $G_L$ exists by Theorem \ref{GL}.  In this section we consider the size and possible values of $g_L$.   We show $g_L\neq G_L$ in general, and that $g_L\geq 7$ for every Gaussian line  $L$. We give necessary and sufficient conditions for $g_L=7$ and describe infinitely many Gaussian lines with $g_L\geq 260{,}000$.  We conjecture that both $g_L$ and $G_L$ can be arbitrarily large. We leave the reader with this and other open problems about $g_L$ and $G_L$.  All of our results are given in terms of which primes are in the divisor set $\DL$ of $L$.  Recall that Theorem~\ref{DL} provides a very simple test for this.

What can we say about $g_L$ for a primitive Gaussian line $L$? In the special case where $L$ is the real or imaginary line, we have 
$g_L=G_L=17$. If $L$ is not the real or imaginary line,  then the  value of $g_L$ will depend on which Gaussian primes occur as divisors of elements on $L$, that is, it will depend on the divisor set $\DL$ defined in Section 2, as well as how the conjugate primes are distributed on the line.
Recall that the divisor set $\DL$ contains only finitely many primes $p\equiv 3\pmod 4$ and for all but finitely many primes $p\equiv 1\pmod 4$, it contains two non-associate conjugate Gaussian primes that lie over $p$ (see Remark 1 following Theorem \ref{DL}).  Intuitively, one expects that $g_L$ will be small if $\DL$ contains lots of small primes, and very large if $\DL$ contains few primes $p\equiv 3\pmod 4$ (or none at all) and exactly one non-associate prime lying over $p$ for small  primes $p\equiv 1\pmod 4$.
We use this idea to show that  that $g_L$ can be as small as $7$, to characterize those Gaussian lines for which $g_L=7$, and to
show that there are infinitely many lines with $g_L\geq 260{,}000$.  As discussed in the introduction, the fact that $g_L\geq 7$ also follows from the main result of Ghorpade and Ram~\cite{ghorpade}.

 \begin{theorem}\label{gLmin}
If $L$ is a primitive Gaussian line, then $g_L\geq 7$.  Moreover, $g_L=7$ if and only if $1+i, 3, 1+2i,1-2i\in\DL$ and the conjugate primes $1+2i$ and $1-2i$ divide consecutive Gaussian integers on $L$.  
\end{theorem}

\begin{proof} Let $L$ be a primitive Gaussian line.  To see that $g_L\geq 7$, consider an arbitrary sequence 
$S=(\alpha_{k+1}, \alpha_{k+2}, \ldots, \alpha_{k+n})$ of $n$ consecutive Gaussian integers on $L$ where $n\leq 6$. 
By Theorem \ref{periodicity}, if two terms of $S$ have a 
common Gaussian prime divisor then that common divisor must lie over a rational prime $p < n$. 
Thus,  we only consider divisibility by $1+i$, $3$, $1+2i$, and $1-2i$ since these are representatives for the non-associate Gaussian primes that lie over 2, 3, and 5.

If $n=2$ or $n=3$, then $\alpha_{k+2}$ is coprime to the other terms of $S$  since two consecutive elements on $L$ are always coprime by Theorem \ref{periodicity}. If $n=4$, then, again by Theorem \ref{periodicity}, there are at exactly two terms of $S$ that are divisible by $1+i$ 
(either $\al_k$ and $\al_{k+2}$, or 
$\al_{k+1}$ and $\al_{k+3}$) and the remaining two terms cannot both be divisible by 3.  Thus, there is at least one term of $S$ that is coprime to  $3(1+i)$, and so is coprime to all  other terms of $S$. If $n=5$, then there are no new Gaussian primes to consider.  Thus, there is similarly  at least one term of $S$  that is coprime to  $3(1+i)$ and so coprime to all other terms of $S$.

If $n=6$, then we must also consider divisibility by the two Gaussian primes that lie over 5. Notice, however, that only $\alpha_{k}$ and $\alpha_{k+5}$ can share a common prime divisor that lies over 5.  Thus, $1+i$ and 3 are the only possible common Gaussian prime divisors of   $\alpha_{k+1}$, 
$\alpha_{k+2}$, $\alpha_{k+3}$, and  $\alpha_{k+4}$ with other terms of $S$. But, as in the case $n=4$, at least one of these four terms is coprime to  $3(1+i)$ and hence is coprime to all other terms of $S$. Thus, for all $n\leq 6$, the sequence $S$ contains at least one term that is coprime to all the others, so $g_L\geq 7$.

Now consider $n=7$.  Suppose that $L$ is a primitive Gaussian line such that $1+i, 3, 1+2i,1-2i\in\DL$ and $1+2i$ and $1-2i$ divide consecutive Gaussian integers on $L$.  Without loss of generality, we assume that if $1+2i$ divides $\alpha_t$, then $1-2i$ divides $\alpha_{t+1}$.  Then by the  Chinese remainder theorem for Gaussian lines (Theorem \ref{crtGL}), there is a $k\in\Z$ such that $(1+i)3(1+2i)$ divides 
$\alpha_{k+1}$. Then, by the periodicity of divisibility (Theorem \ref{periodicity}),  $1-2i$ divides $\alpha_{k+2}$ and $\alpha_{k+7}$, and none of the consecutive Gaussian integers $\alpha_{k+1}, \alpha_{k+2}, \ldots, \alpha_{k+7}$ on $L$ is coprime to all the others.  Thus $g_L=7$.

Conversely, suppose that $L$ is a primitive Gaussian line with  $g_L=7$. Then there is a sequence 
$S=(\alpha_{k+1}, \alpha_{k+2}, \ldots, \alpha_{k+7})$ of seven consecutive Gaussian integers on $L$ where none of the terms  is coprime to all the others. By the periodicity of divisibility, the only possible common Gaussian prime divisors of the terms of $S$ are the associates of  $1+i, 3, 1+2i,1-2i$.
At least two terms of $S$ are not divisible by $1+i$ or $3$, so they must be divisible by primes lying over~5.
They cannot both be divisible by the same prime lying over~5 since they are either 2, 4, or 6 terms apart.  Thus, $1+2i$ and 1$-2i$ must both divide two terms of $S$.  By Theorem \ref{periodicity}, the only possibility is that one divides
 $\alpha_{k+1}$ and $\alpha_{k+6}$ and the other divides 
 $\alpha_{k+2}$ and $\alpha_{k+7}$.  Thus, $1+i, 3, 1+2i,1-2i\in\DL$ and $1+2i$ and $1-2i$ divide consecutive Gaussian integers on $L$.
\end{proof}

\noindent{\em Remark 3.}  It is interesting to consider the problem of characterizing those integers $n$ for which there is a Gaussian line with $g_L=n$. Here we make some remarks on this problem for $n\leq 15$, but the proofs are largely left to the reader because while they are not complicated, they are quite messy.  As a companion to the second half of Theorem~\ref{gLmin}, it is not hard to show that $g_L=9$ if and only if
$1+i, 3, 1+2i,1-2i\in\DL$ and the conjugate primes $1+2i$ and $1-2i$ do {\em not} divide consecutive integers on $L$. Thus, by 
Theorem~\ref{newchinese}, there are infinitely many lines with $g_L=7$ and infinitely many with $g_L=9$.  
If $L$ is a primitive Gaussian line with 
$g_L\neq 7$ or $9$ (so at least one of  $1+i, 3, 1+2i,1-2i$ is not in $\DL$), then an argument similar to that used to prove the first part of Theorem~\ref{gLmin}, shows that any sequence of  $2\leq n\leq 14$ consecutive Gaussian integers on $L$ contains one that is coprime to all the others.  Thus, there are no Gaussian lines with $g_L= 8$, 10, 11, 12, 13, or 14.  But, $g_L=15$ for the infinitely many primitive Gaussian lines $L$ that satisfy $1+i, 3, 1+2i, 7, 11, 2+3i, 2-3i \in\DL$,
$1-2i\not\in\DL$, and the Gaussian primes lying over 13 divide consecutive Gaussian integers on $L$.  Indeed, assuming  that if $2+3i$ divides an integer on $L$  then $2-3i$ divides the next integer on $L$, it follows from the  Chinese remainder theorem for Gaussian lines (Theorem 
\ref{crtGL}) that there is an integer $t$ such that $2, 3, 11, 2+3i$ all divide $\alpha_t$, $2-3i$ divides $\alpha_{t+1}$, $1+2i$ divides $\alpha_{t+2}$, and $7$ divides $\alpha_{t+5}$.  Then, by the periodicity of divisibility (Theorem~\ref{periodicity}), none of the 15 consecutive Gaussian integers
 $\alpha_{t}, \alpha_{t+1}, \alpha_{t+2}, \ldots , \alpha_{t+14}$  on $L$  is coprime to all the others.  
\bigskip

Unlike like the special case of the real line, in general we have $g_L\neq G_L$.  To see  this, we give an example of an infinite family of Gaussian lines that have $g_L=7$ and $G_L>7$.

\bigskip

\noindent {\bf Example 1.} Let $L$ be a primitive Gaussian line with $1+i, 3, 1+2i,1-2i\in\DL$, $7\not\in\DL$, and such that the conjugate primes $1+2i$ and $1-2i$ divide consecutive Gaussian integers on $L$. Then $g_L=7$ by Theorem \ref{gLmin}.  We show that $G_L>7$ by showing that every sequence of eight consecutive Gaussian integers on $L$ contains one that is coprime to all the others.  Suppose for a contradiction that there is a sequence $S$ of eight consecutive Gaussian integers on $L$ where none is coprime to all the others.  By the periodicity of divisibility,  $1+i$ divides every other term of $S$, so there is a  $k\in\Z$ such that the terms of $S$ that are {\em not} divisible by $1+i$ are 
$\alpha_k, \alpha_{k+2}, \alpha_{k+4}, \alpha_{k+6}$. The only possible common Gaussian prime divisors of these terms with another term of $S$ are $3$, $1+2i$, and $1-2i$. Each of $1+2i$ and $1-2i$ divide only one of the terms $\alpha_k, \alpha_{k+2}, \alpha_{k+4}, \alpha_{k+6}$ since they each divide every fifth integer on the line.  Thus 3 must divide two of these terms, which can only happen if 3 divides $\alpha_k$ and $\alpha_{k+6}$.  Thus, $1+2i$ and $1-2i$ must each divide one of $\alpha_{k+2}$ or $\alpha_{k+4}$, and so also divide $\alpha_{k+7}$ and $\alpha_{k-1}$, respectively. Since $S$ consists of 8 terms, including the seven terms $\alpha_k, \alpha_{k+1}, \alpha_{k+2}, \ldots, \alpha_{k+6}$, it must also contain $\alpha_{k+7}$ or $\alpha_{k-1}$, but not both.  Thus, one of $\alpha_{k+2}$ or $\alpha_{k+4}$ is coprime to all the other terms of the sequence and we have reached a contradiction. 
 Notice the requirement that $7\not\in\DL$ is necessary because otherwise there  would be a sequence of eight consecutive Gaussian integers on $L$ where none is coprime to all the others. Namely, by Theorem \ref{crtGL}, there would an integer $t$ such that $1+i, 3, 1+2i$, and $7$ divide $\alpha_t$ and $1-2i$ divides $\alpha_{t+1}$.  Then by the periodicity of divisibility, none of the eight consecutive Gaussian integers $\alpha_t, \alpha_{t+1}, \alpha_{t+2}, \ldots , \alpha_{t+7}$ would be coprime to all the others.

It is interesting to note that if we consider more primes and suppose that $1+i, 3, 1+2i,1-2i, 11, 2+3i, 2-3i\in\DL$, $7\not\in\DL$, and the conjugate primes $1+2i$ and $1-2i$ divide consecutive Gaussian integers on $L$, then $g_L=7$ and $G_L\geq 15$.  We suppress the details, but a similar argument to that given above shows that any sequence of  $8\leq n\leq 14$ consecutive integers on $L$  contains a term that is coprime to all the others.  When $n=15$, however, there is a sequence of $n$ consecutive integers that does not have this property.  Namely, by Theorem \ref{crtGL}, there is an integer $t$ such that $2, 3, 1+2i, 11, 2+3i$ divide $\alpha_t$, $2-3i$ divides $\alpha_{t+1}$, and $1+2i$ divides $\alpha_{t+2}$.  Then by the periodicity of divisibility, 
none of the 15 consecutive Gaussian integers $\alpha_{t}, \alpha_{t+1}, \alpha_{t+2}, \ldots , \alpha_{t+14}$  on $L$ is coprime to all the others. 

\bigskip

Our next example shows that for any bound $B$ there are infinitely many primitive Gaussian lines $L$ such that for all  $7\leq n\leq B$, $L$ contains $n$ consecutive Gaussian integers where none is coprime to all the others.  For these lines, $g_L=7$ and either $G_L=7$ or $G_L>B$.

\bigskip

\noindent {\bf Example 2.}  Let  $B\geq 7$ be an integer. By Theorem \ref{newchinese}, there are infinitely many primitive Gaussian lines $L$ such that 
$1-2i$ divides $\alpha_0$, 
$1+2i$ divides $\alpha_1$, and for all rational primes $p<B$, $p\neq 5$, a Gaussian prime that lies over  $p$ also divides $\alpha_0$.  For these lines, $g_L=7$ by Theorem \ref{gLmin}.  Also, 
$\alpha_1$ and $\alpha_6$ are both divisible by $1+2i$ and, for $2\leq t\leq k-1$,  $\alpha_0$ and $\alpha_t$ are both divisible by a Gaussian prime that divides $t$. Thus, for any integer $n$ with $7\leq n\leq B$, none of the $n$ consecutive Gaussian integers $\alpha_0, \alpha_1, \ldots, \alpha_{n-1}$   on $L$ is  coprime to all the others. Therefore, either $G_L=7$ or $G_L>B$.

\bigskip

Examples 1 and 2  illustrate the intuitive idea that if a primitive Gaussian line $L$ has lots of small primes in its divisor set, then $L$ will contain short sequences of consecutive Gaussian integers where none is coprime to all the others. Conversely,  if the divisor set of $L$ contains few small Gaussian primes, then for small values of $n$, every sequence of $n$ consecutive Gaussian integers on $L$ will contain at least one term that is coprime to all the others.  This idea is made more precise in our final example and subsequent theorem.

\bigskip

\noindent {\bf Example 3.} Let $2\leq k\leq  \numprint{260,000}$.   We show that there are infinitely many Gaussian lines  $L_k$ with $g_{L_k}>k$.  We choose the lines $L_k$ so that their divisor sets contain the minimum possible number of Gaussian primes $\pi$ with $\nu(\pi)\leq k$, where the function $\nu$ is defined in Equation~(\ref{morm}).
Recall from Remark 1 following Theorem \ref{DL}, that the divisor set $\DL$ of any Gaussian line $L$ must contain at least one Gaussian prime lying over $p$ for every rational prime $p\equiv 1\pmod 4$.  We choose $L_k$ so that $\DLk$ contains exactly one Gaussian prime $\pi$ lying over $p$ for each rational prime $p\equiv 1\pmod 4$, $p\leq k$, and 
such that $\DLk$ does not contain $1+i$ or any rational prime $p\equiv 3 \pmod 4$, $p\leq k$.  By Theorem~\ref{newchinese}, there are infinitely many such primitive Gaussian lines $L_k$.  We computationally show that $g_{L_k}\geq k$ for all of these lines. 

Suppose $t$ is a positive integer with $t\leq k\leq \numprint{260,000}$.  Let $S_t$ be a set of $t$ consecutive Gaussian integers on $L_k$. If $\pi$ is a Gaussian prime that divides two distinct elements in $S_t$,  then we must have that $\nu(\pi)< t$  by Theorem \ref{periodicity}.  By the definition of $L_k$, it follows that $N(\pi)=p$ for some rational prime $p\equiv 1\pmod 4$, $p<t$, since these are the only Gaussian primes with $\nu(\pi)< t$ in 
$\DLk$. Now, $\pi$  divides every $p$th Gaussian integer in $S_t$ by Theorem \ref{periodicity}, so it divides at most $\lceil{{t/p}}\rceil$ distinct Gaussian integers in $S_t$, where $\lceil{{\ }}\rceil$ denotes the the ceiling function. This is true for every rational prime  $p\equiv 1\pmod 4$, $p<t$, since $\DLk$ contains a unique Gaussian prime of norm $p$ for each of these prime.  Since these are the only possible Gaussian prime divisors of two elements in $S_t$, it follows that  a crude upper bound for the number of elements of $S_t$ that are divisible by a prime $\pi\in{\cal D}(L_k)$ with $\nu(\pi)< t$ is 
$$B_t=\sum_{{p<t, p\ {\rm prime}}\atop {p\equiv 1\hskip-.1in\pmod 4}}{\left\lceil{t\over p}\right\rceil}.$$
Note that $B_t$ is an overestimate because it over counts the elements in $S_t$ that are divisible by more than one distinct prime $\pi_p\in{\cal D}(L_k)$ with $\nu(\pi_p)=p< t$.
We chose this crude estimate because it is easy to  compute.
For example, if $t=100$, then one can check by hand that $B_{100}=54$.  It follows that if $S_{100}$ is any set of 100 consecutive Gaussian integers on $L_k$, then there are at least 36 elements in $S_{100}$ that are coprime to all the other elements in the set. In particular, $g_{L_k}\neq 100$.  We wrote a program in SageMath \cite{sage} to quickly compute $B_t$ for $t\leq \numprint{260,000}$.  In every case we had $B_t<t$.  
Thus, if $t\leq \numprint{260,000}$, then $S_t$ contains at least one element whose Gaussian prime divisors all satisfy $\nu(\pi)\geq t$.  This element is coprime to all other elements in $S_t$  by Theorem \ref{periodicity}, and so $g_{L_k}\neq t$.  Since this holds for all $t\leq k$,  it follows that $g_{L_k}> k$.

\bigskip

The following theorem is immediate  from Example 3.

\begin{theorem}\label{biggL} There are infinitely many primitive Gaussian lines $L$ for which $g_L\geq \numprint{260,000}$.  
\end{theorem}

The bound $B_t$ in Example 3 actually satisfies $B_t<t$ for $t\leq\numprint{260,185}$, and $B_t>t$ for $t=\numprint{260,186}$. Thus,  there are actually infinitely many primitive Gaussian lines $L$ with $g_L\geq \numprint{260,185}$.  The property that $g_{L_k}\geq k$ could be extended to even larger values of $k$ by lowering $B_t$ to account for elements that are  divisible by more than one Gaussian prime lying over a rational prime $p\equiv 1\pmod 4$. It is doubtful, however, that this method would lead to a proof that $g_{L_k}\geq k$ for all $k$.  Yet it seems plausible that this is true and, in particular, that $g_L$ can be arbitrarily large.  We end with this and other open problems.

\bigskip

\noindent{\bf Open problems.} The following are among the problems about $G_L$ and $g_L$ that arise from our work:
\begin{enumerate}
\item Find the value of $G_L$ for some primitive Gaussian line with $\Delta\neq 0$.
\item Characterize those integers $n$ for which there is a primitive Gaussian line with $g_L=n$.
In Remark 3, we note that there are infinitely many Gaussian lines with  $g_L=7, 9$ and $15$, and none with $g_L=8, 10, 11, 12, 13, 14$.  What other values of $g_L$ are possible?  For instance, can $g_L$ be even? 
\item Determine whether or not $G_L$ can be arbitrarily large. For every $n\in\N$ is there a primitive Gaussian line $L$ with $G_L\geq n$?
\item Determine whether or not $g_L$ can be arbitrarily large. For every $n\in\N$ is there a primitive Gaussian line $L$ with $g_L\geq n$?

\end{enumerate}
We suspect that $g_L$, and hence $G_L$, can be arbitrarily large (for instance, see Example 3), but perhaps only $G_L$ can be arbitrarily large and there is an upper bound $B\in\Z$ such that $g_L\leq B$ for all Gaussian lines $L$.  It is also possible that both $g_L$ and $G_L$ are bounded. It would be interesting to know if they are bounded, and moreover, what integer values of $g_L$ and of $G_L$ are possible. Determining the possible values seems considerably harder for $G_L$ than for $g_L$ since we only know its value for the real and imaginary lines.
We leave the reader with these challenges.

\end{document}